\documentclass{commat}

\newcommand{\norm}[1]{\ensuremath{\left\|#1\right\|}}

\title{%
    Polynomial complex Ginzburg-Landau equations in almost periodic spaces
}

\author{%
    Agustín Besteiro
}

\affiliation{\address{Universidad Abierta Interamericana, Centro de Altos Estudios en Tecnología Informática, Buenos Aires, Argentina
}
    \email{%
    agustin.besteiro@uai.edu.ar
}
}

\abstract{%
We consider complex Ginzburg-Landau equations with a polynomial nonlinearity in the real line. We use splitting-methods to prove well-posedness for a subset of almost periodic spaces. Specifically, we prove that if the initial condition has multiples of an irrational phase, then the solution of the equation maintains those same phases.
}



\VOLUME{31}
\YEAR{2023}
\NUMBER{1}
\firstpage{91}
\DOI{https://doi.org/10.46298/cm.10279}

\begin{paper}

\section{Introduction}
We consider the 1-dimensional autonomous system
\begin{equation}\label{eq: CGL}
\begin{cases}
\partial_t u = (\alpha +i \beta) \partial_{xx} u + \gamma u + (a+ib)B(u), \\
u(0) = u_0,
\end{cases}
\end{equation}
where $u(x,t)$ is a complex valued function with $x\in\mathbb{R}$, $t>0$, $\alpha>0$, $\beta>0$, $\gamma\geq 0$, $a>0$, $b>0$ and $B$ a continuous map.
The linear term represented by $(\alpha +i \beta) \partial_{xx}$ characterizes the Complex Ginzburg-Landau equations.
For $\beta = 0$~\eqref{eq: CGL} reduces to a nonlinear heat equation and for $\alpha = 0$ to a nonlinear Schrödinger equation.
A large amount of work has been done to prove well-posedness of~\eqref{eq: CGL} with different nonlinearities (see, for instance,~\cite{Aranson2002},~\cite{Ginibre1996},~\cite{Ginibre1997}).
In our case, we study well-posedness of problem~\eqref{eq: CGL} with polynomial nonlinearities. These nonlinearities are considered in Fisher-Kolmogorov equations and Fitzhugh-Nagumo equations~\cite{Asgari2011},~\cite{Fisher1937},~\cite{Kolmogorov1989}. Almost periodic spaces were introduced by Bohr~\cite{Bohr2018} and further developed by Stepanoff~\cite{Stepanoff1925}, Weyl~\cite{Weyl1927} and Besicovitch~\cite{Besicovitch1926},~\cite{Besicovitch1954}. These spaces are well-studied for the Complex Ginzburg-Landau equations with different nonlinearities~\cite{Guo2000},~\cite{Gao2018}.
We consider $u_0$ with specific irrational phases, and we prove that the time evolution by~\eqref{eq: CGL} maintains the same phases. We use splitting-methods for evolution equations developed for numerical purposes~\cite{DeLeo2015},~\cite{Borgna2015}. These methods were used to prove well-posedness of Complex Ginzburg-Landau equations and Reaction-diffusion equations in other spaces~\cite{Besteiro2019},~\cite{Besteiro2020},~\cite{Besteiro2018},~\cite{besteiro1}.

The paper is organized as follows: In Section 2 we set notations and preliminary results. In Section 3 we analyze the nonlinear problem. Finally, in Section 4 using splitting methods, we combine results from Sections 2 and 3 to obtain that the solution of~\eqref{eq: CGL} is in a subset of an almost periodic space.

\section{Notations and Preliminaries.}\label{sec:2}
We introduce some definitions and preliminary results.
\begin{definition}
We define $C_\textrm{u}(\mathbb{R})$ as
the set of uniformly continuous and bounded functions on $\mathbb{R}$ equipped with the norm,
\begin{align*}
    \norm{u}_\infty = \sup_{x\in \mathbb{R}}|u(x)|.
\end{align*}
\end{definition}

\begin{definition}
We define the set of almost periodic functions~\cite{Bohr2018} as,
\begin{align*}
    P(\mathbb{R}) = \{u\in C_u(\mathbb{R}): \ u = \sum^{\infty}_{j = 1}a_j e^{ix\lambda_j} :\lambda_j \in \mathbb{R} \}
\end{align*}
equipped with the uniform continuous norm.
\end{definition}

\begin{definition}
We define the following subset of almost periodic functions,
\begin{align}\label{def: almostperiodic}
    A_\lambda(\mathbb{R}) = \{u\in C_u(\mathbb{R}): \ u = \sum^{\infty}_{j = 1}a_j e^{ixj\lambda};\sum^{\infty}_{j = 1}|a_j|<\infty \ : \lambda
    \in \mathbb{R} \}
\end{align}
equipped with the norm:
\begin{align}\label{def: normalmostperiodic}
\norm{u}_{A_{\lambda}} = \sum^{\infty}_{j = 1}|a_j|.
\end{align}
\end{definition}
\begin{theorem}
 $A_{\lambda}(\mathbb{R})$ is a Banach space.
\end{theorem}
\begin{proof}
    Let $\{u_n\}_{n\in \mathbb{N}}$ a Cauchy sequence such that $u_n \in A_{\lambda}(\mathbb{R}) \subset P(\mathbb{R})$. Then $\{u_n\}_{n\in \mathbb{N}}$ is a bounded sequence, that is, $\norm{u_n}_{A_{\lambda}} = \sum^{\infty}_{j = 1}|a^n_j|\leq K$ for all $n \in \mathbb{N}$. Additionally, $\{u_n\}_{n\in \mathbb{N}}$ is a convergent sequence in $P(\mathbb{R})$, i.e. $u_n \to u$ with $u \in P(\mathbb{R})$ and $\norm{u}_{A_{\lambda}} = \sum^{\infty}_{j = 1}|a_j|$. Therefore, we have that $\sum^N_{j = 1}|a_j| = \lim_{n \to \infty} \sum^N_{j = 1}|a^n_j| \leq K$, which implies that 
    \[
    \lim_{N \to \infty} \sum^N_{j = 1}|a_j| = \sum^{\infty}_{j = 1}|a_j| \leq K.
    \]

    On the other hand, we can define an inner product for $u,v \in P(\mathbb{R}$):
    \begin{align*}
        \langle u(x), v(x) \rangle = \lim_{T\to \infty}\frac{1}{T}\int_{0}^{T}u(x)v(x)dx.
    \end{align*}
    In particular, we have a normalized orthogonal system in the following sense (see~\cite{Bohr2018}):
    \begin{align*}
        \langle e^{i\lambda_{1}x}, e^{-i\lambda_{2}x} \rangle
        = \lim_{T \to \infty}\frac{1}{T}\int_{0}^{T}e^{i\lambda_{1}x}e^{-i\lambda_{2}x}
        =
\begin{cases}
        0 & \text{, for } \lambda_1 \neq \lambda_2 \\
        1 & \text{, for } \lambda_1 = \lambda_2.
        \end{cases}
    \end{align*}
    For $u \in P(\mathbb{R}),\ u = \sum_{j = 1}^{\infty} a_j e^{i \lambda_j x}$, we have,
    \begin{align*}
        \langle u(x), e^{-i\lambda_jx} \rangle = a_j.
    \end{align*}
    If $u_n \in A_{\lambda}(\mathbb{R})$ is a Cauchy sequence, then
    \begin{align*}
        a_{j}^{n}
        = \langle u_n(x), e^{-ix j \lambda}\rangle \ \xrightarrow{n \to \infty} \ \langle u(x),e^{-i x j\lambda} \rangle 
        = a_{j}.
    \end{align*}
    As~\eqref{def: almostperiodic} and~\eqref{def: normalmostperiodic} are met then $u \in A_{\lambda}(\mathbb{R})$.
\end{proof}

\begin{remark}
As $\ell ^1 (\mathbb{N})$ is a Banach Algebra then, if $u \in A_{\lambda}(\mathbb{R})$ and $v \in A_{\lambda}(\mathbb{R})$ we have that,
\begin{align*}
    \norm{uv}_{A_\lambda}\leq \norm{u}_{A_\lambda}\norm{v}_{A_\lambda}.
\end{align*}
\end{remark}
The following definitions and proofs can be extended to $x \in \mathbb{R}^d$ (See~\cite{Engel2001}).

\begin{definition}
We denote $U(t)$ as the one parameter semigroup that solves the underlying linear equation
\begin{align}
    \partial_t u = (\alpha + i \beta) \partial_{xx}u + \gamma u.
\end{align}
The operator can be represented by the convolution in x
\begin{align*}
    U(t) = {(4\pi t(\alpha +i \beta))}^{-1/2} e^{(-x^2 /[4t(\alpha+i\beta)]) + \gamma t}*u_0 = G_t(x)*u_0
\end{align*}
and the kernel $G_t$ satisfies:
\begin{align*}
    |G_t(x)| = (4\pi t{(\alpha^2 +\beta^2)}{(4\pi t(\alpha^2 +\beta^2)^{1/2})}^{-1/2} e^{(-x^2 /[4t(\alpha^2 +\beta^2)]) + \gamma t}.
\end{align*}
Clearly, $G_t(x) \in L^1 (\mathbb{R})$.
\end{definition}

\begin{proposition}
For each $t \ge 0$, define $U(t)u_0 = G_t*u_0$.  The one-parameter family of operators $\{U(t)\}_{t\geq0}$ ia a strongly continuous semigroup in $C_u(\mathbb{R})$.
\end{proposition}
\begin{proof}
The semigroup property, $U(t)U(t')u = U(t+t')u$ is proven similarly to the heat kernel. We show that, $U(t)u$ converges to $u$ for all $u\in C_\textrm{u}(\mathbb{R})$ when $t\to 0$. Indeed, we have,
\begin{align*}
|(U(t)u)(x)-u(x)|
\le{ }&{ }\int_{\mathbb{R}}G_{t}(y)|u(x-y)-u(x)|dy \\
={ }&{ } \int_{|y|<\delta}G_{t}(y)|u(x-y)-u(x)|dy
+\int_{|y|\ge \delta}G_{t}(y)|u(x-y)-u(x)|dy.
\end{align*}
The first integral of the right side of the equality can be estimated as follows:
\begin{align*}
\int_{|y|<\delta}G_{t}(y)|u(x-y)-u(x)|dy
& \le \int_{\mathbb{R}}G_{t}(y)\max_{|y|<\delta}|u(x-y)-u(x)|dy \\
& = \max_{|y|<\delta}|u(x-y)-u(x)|.
\end{align*}
This can be small enough because, $|y|<\delta$ and $u$ is uniformly continuous. For the second term we proceed in the following way,
\begin{align*}
\int_{|y|\ge \delta}G_{t}(y)|u(x-y)-u(x)|dy
& \leq 2\|u\|_{\infty} \int_{|y|\ge \delta} G_{t}(y) dy.
\end{align*}
Since $(-x^2 /[4t(\alpha^2 +\beta^2)]) \to\infty$ when $t\to 0^{+}$ and $G_{t}\in L^{1}(\mathbb{R})$, the right side of the previous equality tends to $0$.
The next property proves that $U$ is well defined, that is $U(t)u\in C_\textrm{u}(\mathbb{R})$.
\begin{align*}
|(U(t)u)(x_{1})-(U(t)u)(x_{2})| & \le
\int_{\mathbb{R}}G_{t}(y)|u(x_{1}-y)-u(x_{2}-y)|dy \\
& \le \varepsilon\int_{\mathbb{R}}G_{t}(y)dy = \varepsilon,
\end{align*}
In the last inequality we used that $u$ is uniformly continuous.
\end{proof}

\begin{lemma}\label{le: linearpol}
If $u_0 \in A_{\lambda}(\mathbb{R})$ then $U(t)u_0 \in A_{\lambda}(\mathbb{R})$ for $t>0$
\end{lemma}
\begin{proof}
    As $u_0 \in A(\mathbb{R})$ then we have
    \begin{align*}
       U(t)u_0 = G_t*u_0 =
       \int_{\mathbb{R}}{(4\pi t(\alpha +i \beta))}^{-1/2} e^{(-y^2 (\alpha-i\beta)/[4t(\alpha^2 +\beta^2)]) + \gamma t}u_0(x-y)dy
    \end{align*}
   where
   \begin{align}
       u_0(x-y) = \sum_{j = 1}^{\infty}a_j e^{ij \lambda(x-y)}
   \end{align}
   then we have,
   \begin{align*}
       U(t)u_0 
       { }&{ }= \int_{\mathbb{R}}{(4\pi t(\alpha +i \beta))}^{-1/2} e^{(-y^2 (\alpha-i\beta)/[4t(\alpha^2 +\beta^2)]) + \gamma t}\sum_{j = 1}^{\infty}a_j e^{ij \lambda(x-y)}dy \\
       { }&{ }= \int_{\mathbb{R}}{(4\pi t(\alpha +i \beta))}^{-1/2}\sum_{j = 1}^{\infty}a_j e^{(-y^2 (\alpha-i\beta)/[4t(\alpha^2 +\beta^2)]) + \gamma t +ij \lambda(x-y)}dy \\
       { }&{ }= {(4\pi t(\alpha +i \beta))}^{-1/2}\int_{\mathbb{R}}\sum_{j = 1}^{\infty}a_j B(y) e^{ij \lambda{(x-y)}}dy,
   \end{align*}
   with
   \begin{align*}
       B(y) = e^{(-y^2 (\alpha-i\beta)/[4t(\alpha^2 +\beta^2)]) + \gamma t}.
   \end{align*}
   Using dominated convergence theorem and that
   \begin{align*}
   \left|\sum_{j = 1}^n a_j B(y) e^{ij \lambda{(x-y)}}\right|
   \leq \sum_{j = 1}^n |a_j B(y)|
   = |B(y)|\sum_{j = 1}^n |a_j|\leq K |B(y)|
   \end{align*}
   we have,
 \begin{align*}
           U(t)u_0 = {(4\pi t(\alpha +i \beta))}^{-1/2}\sum_{j = 1}^{\infty}\int_{\mathbb{R}}a_j B(y) e^{ij \lambda{(x-y)}}dy.
 \end{align*}

 Finally, we know that
 \begin{align*}
     \int_{\mathbb{R}}e^{-ax^2 -bx+c}dx = \frac{e^{b^2 /(4 a)+c}\sqrt{\pi}}{\sqrt{a}}
 \end{align*}
 with $Re(a)>0$. In our case we have
 \[
 a = (\alpha-i\beta)/[4t(\alpha^2 +\beta^2)]) = 1/(4t(\alpha+i\beta)),
\quad
 b = ij \lambda
 \quad \textup{ and } \quad
 c = ixj \lambda+\gamma t.
 \]
       Then we have,
       \begin{align*}
           U(t)u_0 =& {(4\pi t(\alpha +i \beta))}^{-1/2}\sum_{j = 1}^{\infty}a_j\frac{e^{-t{(j \lambda)}^2 (\alpha + i \beta)+ixj \lambda+\gamma t}\sqrt{\pi}}{\sqrt{(1/(4t(\alpha+i\beta))}}
           \\
           = & \sum_{j = 1}^{\infty}a_je^{-t{(j \lambda)}^2 (\alpha + i \beta)+ixj \lambda+\gamma t} \\
           = &\sum_{j = 1}^{\infty}A_je^{ixj \lambda},
       \end{align*}
      where $A_j = a_je^{-t{(j \lambda)}^2 (\alpha + i \beta)+\gamma t}$, so we have~\eqref{def: almostperiodic}. On the other hand we have,
      \begin{align*}
      \sum_{j = 1}^{\infty}|A_j| =
      \sum_{j = 1}^{\infty}|a_je^{-t{(j \lambda)}^2 (\alpha + i \beta)+\gamma t}| =
      \sum_{j = 1}^{\infty}|a_je^{-t{(j \lambda)}^2 \alpha+\gamma t}|<\infty
      \end{align*}
      Then we have~\eqref{def: normalmostperiodic} and $U(t)u_0 \in A_{\lambda}(\mathbb{R})$
\end{proof}

Next, we consider integral solutions of the problem~\eqref{eq: CGL}.

We say that $u\in C([0,T],C_\textrm{u}(\mathbb{R}))$ is a mild solution of~\eqref{eq: CGL}
if and only if $u$ verifies
\begin{align}\label{eq: mild solution}
u(t) = U(t)u_{0} + \int_{0}^{t} U(t-t')B(u(t')) dt'.
\end{align}

If $F$ is a locally Lipschitz map, for any $z_{0} \in C_\textrm{u}(\mathbb{R})$
there exists a unique solution of the equation
\begin{align}\label{eq: ode}
\begin{cases}
 \partial_tz = B(z), \\
 z(0) = z_{0},
 \end{cases}
\end{align}
defined in the interval $[0,T^{*}(z_{0}))$.
Moreover, there exists a function $\bar{T}:[0,\infty)\to [0,\infty)$, which is non-increasing and such that $T^{*}(z_{0}) \ge \bar{T}(|z_{0}|)$.
The solution of~\eqref{eq: ode} is solution of the integral equation
\begin{align}\label{eq: integral equation}
 z(t) = z_{0} + \int_{0}^{t} B(z(t')) dt'.
\end{align}
Also, one of the following alternatives holds:
\begin{itemize}
\item[-] $T^{*}(z_{0}) = \infty$;
\item[-] $T^{*}(z_{0}) < \infty$ and $|z(t)| \to \infty$ when $t \uparrow T^{*}(z_{0})$.
\end{itemize}

We will denote by $\mathsf{N}(t,.):C_\textrm{u}(\mathbb{R})\to C_\textrm{u}(\mathbb{R})$ the flow
generated by the ordinary equation, i.e.: for any $x\in\mathbb{R}$,
$\mathsf{N}(t,u_{0})(x)$ is the solution of the problem~\eqref{eq: ode}
with initial data $z_{0} = u_{0}(x)$. Therefore, if $u(t) = \mathsf{N}(t,u_{0})$
\begin{align*}
 u(x,t) = u_{0}(x) + \int_{0}^{t}B(u(x,t'))dt'.
\end{align*}

We recall well-known local existence results for evolution equations.

\begin{theorem}\label{th: local existence}
	There exists a function $T^{*}:C_\textrm{u}(\mathbb{R})\to \mathbb{R}_{+}$ such that
	for $u_{0}\in C_\textrm{u}(\mathbb{R})$, exists a unique $u\in C([0,T^{*}(u_{0})),C_\textrm{u}(\mathbb{R}))$ mild solution
	of~\eqref{eq: CGL} with $u(0) = u_{0}$. Moreover, one of the following alternatives holds:
	
\begin{itemize}
		\item $T^{*}(u_{0}) = \infty$;
		\item $T^{*}(u_{0}) < \infty$ and $\lim_{t \uparrow T^{*}(u_{0})}|u(t)| = \infty$.
	
\end{itemize}
\end{theorem}
\begin{proof}
	See Theorem 4.3.4 in~\cite{Cazenave1998}.
\end{proof}

\begin{proposition}\label{pr: continuous dependence}
	Under conditions of the theorem above, we have the following statements:
	
\begin{enumerate}
		\item $T^{*}:C_\textrm{u}(\mathbb{R})\to \mathbb{R}_{+}$ is lower semi-continuous;
		\item If $u_{0,n} \to u_{0}$ in $C_\textrm{u}(\mathbb{R})$ and $0 < T < T^{*}(u_{0})$, then
		$u_{n} \to u$ in the Banach space $C([0,T],C_\textrm{u}(\mathbb{R}))$.
	
\end{enumerate}
\end{proposition}
\begin{proof}
	See Proposition 4.3.7 in~\cite{Cazenave1998}.
\end{proof}

\section{Nonlinear equation}\label{sec:3}
In order to apply the Lie-Trotter method, we prove that if the initial state $u_0 \in A_{\lambda}(\mathbb{R})$ the solution of the nonlinear ordinary equation $z(t) \in A_{\lambda}(\mathbb{R})$. We consider the equation with a cubic nonlinearity and then we extend the result to a nth-degree nonlinearity.

We study first, the solution for the nonlinear equation~\eqref{eq: ode} with a cubic nonlinearity, that is
\begin{align}\label{eq: ode2}
\begin{cases}
 \partial_tz = -(a+ib)z^3, \\
 z(0) = z_{0}.
 \end{cases}
\end{align}
\begin{lemma}\label{le: nonlinearcubic}
    If $u_0(x) = z_0 \in A_{\lambda}(\mathbb{R})$ then the solution of the equation~\eqref{eq: ode2}, $z(t) \in A_{\lambda}(\mathbb{R})$ for $t \in (0,T^* (z_0))$.
\end{lemma}
  \begin{proof}
  We prove that $F:A_{\lambda}(\mathbb{R}) \to A_{\lambda}(\mathbb{R})$ and that $F$ is a locally Lipschitz map in $A_{\lambda}(\mathbb{R})$.
  Let $u \in A_{\lambda}(\mathbb{R})$ and $k = -(a+ib)$ then we have,
  \begin{align*}
      F(u) = ku^3 = \ & k \sum_{j = 1}^{\infty}\sum_{k = 1}^{\infty}\sum_{l = 1}^{\infty}a_ja_ka_l e^{ix(j+k+l)\lambda} \\
= & \ k \sum_{j = 1}^{\infty}\sum_{m = 1}^{\infty}\sum_{l = 1}^{\infty}a_ja_{m-j-l}a_l e^{ixm\lambda} = \sum_{m = 1}^{\infty}A_m e^{ixm\lambda}
  \end{align*}
  with $m = j+k+l$ and $A_m = {(a_{m})}^3$. Also we have that $\sum_{m = 1}^{\infty}|A_m| = \sum_{m = 1}^{\infty}|a_ja_ka_l|<\infty$.
  On the other hand, if $u,v\in A_{\lambda}(\mathbb{R})$ we can see that,
  \begin{align*}
      \norm{F(u)-F(v)}_{A_{\lambda}}
      { }&{ }= \norm{u^3 -v^3}_{A_{\lambda}} \\
      { }&{ }\leq \frac{1}{2} \norm{(u^2 -v^2)}_{A_{\lambda}}\norm{(u+v)}_{A_{\lambda}}+\norm{(u^2 +v^2)}_{A_{\lambda}}\norm{(u-v)}_{A_{\lambda}}.
  \end{align*}
   We use that, $\norm{u+v}_{A_{\lambda}(\mathbb{R})}\leq\norm{u}_{A_{\lambda}(\mathbb{R})}+\norm{v}_{A_{\lambda}(\mathbb{R})} = \sum_{j = 1}^{\infty}|a_j|+\sum_{k = 1}^{\infty}|a_k|<\infty$, a similar procedure proves that $\norm{u^2 +v^2}_{A_{\lambda}}<\infty$ and $\norm{u^2 -v^2}_{A_{\lambda}}<\infty$. Then,
\[
\norm{F(u)-F(v)}_{A_{\lambda}} \leq \frac{1}{2}C \norm{u-v}_{A_{\lambda}}.
\qedhere
\]
\end{proof}

  We generalize the ODE with an n-th degree nonlinearity for $n>2$.

\begin{align}\label{eq: odenth}
\begin{cases}
 \partial_tz = -(a+ib)z^n, \\
 z(0) = z_{0}.
 \end{cases}
\end{align}
\begin{lemma}\label{le: nonlinearnth}
    If $u_0(x) = z_0 \in A_{\lambda}(\mathbb{R})$ then the solution of the equation~\eqref{eq: odenth}, $z(t) \in A_{\lambda}(\mathbb{R})$ for $t \in (0,T^* (z_0))$.
\end{lemma}
  \begin{proof}
 The proof is similar to the previous proof, using that,
 \begin{align*}
     u^n = \sum_{j_1 = 1}^{\infty}\sum_{j_2 = 1}^{\infty}\cdots\sum_{j_n = 1}^{\infty}a_{j_1}a_{j_2} \cdots a_{j_n} e^{ix \lambda\sum_{i = 1}^{\infty}j_i}
 \end{align*}
 and
 \[
     a^n -b^n = (a-b)\left(a^{n-1}+a^{n-2}b+a^{n-3}b^2 +\dots+a^2 b^{n-3}+ab^{n-2}+b^{n-1}\right).
     \qedhere
\]
\end{proof}

\section{Splitting method}
This section is based on the splitting method developed in~\cite{DeLeo2015}. We apply the Lie-Trotter method to the linear and nonlinear problem.
The temporal variable must be broken down into regular intervals and the evolution of the linear and nonlinear problems are considered alternately. This is described by two sequences $\{V_{h,k}\}$ for the linear equation and $\{W_{h,k}\}$ for the nonlinear equation. Using Theorem 3.9 from~\cite{DeLeo2015}, this approximate solution converges to the solution of problem~\eqref{eq: CGL}, when the time intervals $h = t/n \to 0$.

	Let $X$ be a Banach space and we define $\alpha:\mathbb{R} \to \mathbb{R}$ a periodic function of period $1$ as:
\begin{align*}
	\alpha(t) =
	\left\{
\begin{array}{cl}
	2 & \text{, if } k \le t < k+1/2, \\
	0 & \text{, if } k-1/2 \le t < k,
\end{array}
	\right.
\end{align*}
	for $k\in \mathbb{Z}$.

	\label{def: alpha_h}
	Given $h>0$, we define the function $\alpha_{h}:\mathbb{R} \to \mathbb{R}$ as $\alpha_{h}(t) = \alpha(t/h)$.
	Clearly $0\le \alpha_{h} \le 2$, $\alpha_{h}$ is $h$-periodic and its mean value is $1$.

	We consider $\tau_{h}:\mathbb{R}^2 \to \mathbb{R}$ given by
\begin{align*}
	\tau_{h}(t,t') = \int_{t'}^{t}\alpha_{h}(t'')dt''.
\end{align*}
	We define $\Omega = \{(t,t') \in \mathbb{R}^2 : 0 \le t' \le t \}$ and $U_{h}:\Omega \to \mathcal{B}(X)$ given by
	$U_{h}(t,t') = U(\tau_{h}(t,t'))$.

	We consider the system,
\begin{align*}
\begin{cases}
	\partial_t u_h + \alpha_h(t) \partial_{xx} u_h(x,t) = (2-\alpha_h(t)) F(u_h(x,t)), \\
	u_h(x,0) = u_{h0}(x),
\end{cases}
\end{align*}
	where $u(x,t)\in X, \ t>0$ and $\ F:X\to X$ is a continuous function.

	Similarly, we define the integral equation:
\begin{align}\label{eq: u aprox}
	u_{h}(t) & = U_{h}(t,0)u_{h0} + \int_{0}^{t} (2-\alpha_{h}(t')) U_{h}(t,t') F(u_{h}(t')) dt'.
\end{align}
	The following two theorems are a consequence of sections 2 and 3 of~\cite{DeLeo2015}.
\begin{theorem}\label{th: pasospropagador}
	Let $u_{h}$ the solution of~\eqref{eq: u aprox}, if $W_{h,k} = u_{h}(kh)$ y $V_{h,k} = u_{h}(kh-h/2)$,
	then
\begin{subequations}
\begin{align}\label{eq: metodo LT lineal}
		V_{h,k+1} & = U(h)U_{h,k}, \\
		\label{eq: metodo LT no lineal}
		W_{h,k+1} & = N(kh+h,kh+h/2,V_{h,k+1}),
\end{align}
\end{subequations}
	where $N$ is the flux associated to $2F$, that is:
\begin{align*}
\begin{cases}
	\dot{w} = 2F(w(t)), \\
	w(0) = w_{0}.
\end{cases}
\end{align*}
\end{theorem}

\begin{proof}
	For $t_{1}\in (0,t)$ it verifies
\begin{align*}
	u_{h}(t)
	& = U_{h}(t,t_{1}) u_{h0}(t_{1}) + \int_{t_{1}}^{t} (2-\alpha_{h}(t')) U_{h}(t,t') F(u_{h}(t')) dt'
\end{align*}
	using that $t_{1} = kh$ y $t = kh+h/2$, we have
\begin{align*}
	V_{h,k+1} & = U_{h}(kh+h/2,kh) W_{h,k}
	+ \int_{kh}^{kh+h/2} (2-\alpha_{h}(t')) U_{h}(kh+h/2,t') F(u_{h}(t')) dt',
\end{align*}
	given that $\alpha_{h}(t) = 2$ for $t\in [kh,kh+h/2)$, we have $\tau_{h}(kh+h/2,kh) = h$ and therefore~\eqref{eq: metodo LT lineal}.
	Similarly, $\alpha_{h}(t) = 0$ for $t\in [kh+h/2,kh+h)$, then $\tau_{h}(t,kh+h/2) = 0$
	and therefore
\begin{align*}
	u_{h}(t) = V_{h,k+1} + 2 \int_{kh+h/2}^{t}F(u_{h}(t')) dt',
\end{align*}
	evaluating in $t = kh+h$, we obtain~\eqref{eq: metodo LT no lineal}.
\end{proof}
\begin{theorem}\label{th: convergencia}
	Let $u\in C([0,T^{*}),X)$ the solution of the integral problem~\eqref{eq: mild solution}
\begin{align*}
	u(t) = U(t)u_0 + \int_{0}^{t}U(t-t')F(u(t')) dt',
\end{align*}
	$T\in(0,T^{*})$ and $\varepsilon>0$. There exists $h^{*}>0$ such that if $0<h<h^{*}$, then $u_{h}$ the solution of~\eqref{eq: u aprox} with $u_h(x,0) = u_{0}(x)$, is defined in the interval $[0,T]$ and verifies
 \[
 \Vert u(t)-u_{h}(t)\Vert_{X} \le \varepsilon
 \quad \textup{ for } \ t\in [0,T].
 \]
\end{theorem}

\begin{proof}
    See Theorem 3.9 from~\cite{DeLeo2015}.
\end{proof}

We now apply Lemma~\ref{le: linearpol} from Section~\ref{sec:2} related to linear equation and Lemmas~\ref{le: nonlinearcubic} and~\ref{le: nonlinearnth} from Section~\ref{sec:3} related to the nonlinear equations. In order to obtain well-posedness results for the solution $u(t)$ of equation~\eqref{eq: CGL} in $A_{\lambda}(\mathbb{R})$, we use Theorem~\ref{th: convergencia} to join the linear and nonlinear results. The following theorem is proved for the cubic case but the other cases are similar.

\begin{theorem}
 Let $u_{0} \in A_{\lambda}(\mathbb{R})$, then the solution of~\eqref{eq: CGL} $u(t) \in A_{\lambda}(\mathbb{R})$ for $t\in (0,T^* (u_0))$.
\end{theorem}

\begin{proof}
For $t \in [0,T^{*}(u_{0}))$, let $n\in\mathbb{N}$, $h = t/n$ and $\{W_{h,k}\}_{0\le k\le n},\{V_{h,k}\}_{1\le k\le n}$ be the sequences given by $W_{h,0} = u_{0}$,
\begin{subequations}\label{eq: L-T}
\begin{align}
V_{h,k+1} & = U(h) W_{h,k}, \\
W_{h,k+1} & = \mathsf{N}(h,V_{h,k+1}), \quad k = 0,\dots,n-1.
\end{align}
\end{subequations}

     We claim that $W_{h,k+1}\in A_{\lambda}(\mathbb{R})$ for $k = 0,\dots,n$.
 Clearly, the assertion is true for $k = 0$.
 If $W_{h,k}\in A_{\lambda}(\mathbb{R})$, from Lemma~\ref{le: linearpol},
 we have $U(h) W_{h,k} \in A_{\lambda}(\mathbb{R})$. Using Lemma~\ref{le: nonlinearcubic}, we can see that
\begin{align*}
 W_{h,k+1}
= \mathsf{N}(h,V_{h,k})\in A_{\lambda}(\mathbb{R}).
\end{align*}
By theorem~\ref{th: convergencia} we have that
$W_{h,n}\to u(t)$ when $n \to \infty$.

As $A_{\lambda}(\mathbb{R})$ is a Banach space,
we obtain the result.
\end{proof}

\section*{Acknowledgments}
This work was supported by Universidad Abierta Interamericana (UAI), CONICET--Argentina and the project PICTO-2021-UNGS-00001.

\EditInfo{November 15, 2020}{April 30, 2021}{Serena Dipierro}

\end{paper}
\begin{references}

\refer{Paper}{Aranson2002}
\Rauthor{Aranson I. S. and Kramer L.}
\Rtitle{Class of cellular automata for reaction-diffusion systems}
\Rjournal{Rev. Mod.
Phys.}
\Rvolume{74}
\Ryear{2002}
\Rnumber{99}
\Rpages{1749-1752}

\refer{Paper}{Asgari2011}
\Rauthor{Asgari Z. and Ghaemi M. and Mahjani M. G.}
\Rtitle{Pattern Formation of the FitzHugh-Nagumo Model:Cellular Automata Approach}
\Rjournal{Iran. J. Chem. Chem. Eng.}
\Rvolume{30}
\Ryear{2011}
\Rnumber{1}
\Rpages{135-142}

	\refer{Book}{Besicovitch1954}
    \Rauthor{Besicovitch A.S.}
  \Rtitle{Almost periodic functions}
  \Rvolume{4}
  \Ryear{1954}
  \Rpublisher{Dover New York}

\refer{Paper}{Besicovitch1926}
    \Rauthor{Besicovitch A. S.}
  \Rtitle{On generalized almost periodic functions}
   \Rjournal{Proceedings of the London Mathematical Society}
  \Rvolume{2}
  \Ryear{1926}
  \Rnumber{1}
  \Rpages{495-512}

 	\refer{Paper}{besteiro1}
 	\Rauthor{Besteiro A. and Rial D.}
  \Rtitle{Existence of Peregrine type solutions in fractional reaction-diffusion equations.}
  \Rjournal{Electron. J. Qual. Theory Differ. Equ.}
  \Rvolume{9}
	\Ryear{2019}
	\Rpages{1-9}

	\refer{Paper}{Besteiro2018}
  \Rtitle{Global existence for vector valued fractional reaction-diffusion equations}
  \Rauthor{Besteiro A. and Rial D.}
  \Rjournal{Publicacions Matemàtiques}
  \Rpages{653-680}
  \Rvolume{96}
  \Rnumber{2}
  \Ryear{2021}

\refer{Paper}{Besteiro2019}
  \Rtitle{Polynomial complex Ginzburg-Landau equations in Zhidkov spaces}
  \Rauthor{Besteiro A.}
  \Rjournal{Matematychni Studii}
  \Rvolume{52}
  \Rnumber{1}
  \Rpages{55-62}
  \Ryear{2019}

\refer{Paper}{Besteiro2020}
  \Rtitle{A note on dark solitons in nonlinear complex Ginzburg-Landau equations}
  \Rauthor{Besteiro A.}
  \Rjournal{Mathematica}
  \Rvolume{62-85}
  \Rnumber{1}
  \Rpages{11-15}
  \Ryear{2020}

\refer{Book}{Bohr2018}
\Rauthor{Bohr H.}
\Rtitle{Almost periodic functions}
\Rpublisher{Courier Dover Publications}
\Ryear{2018}
\Rpages{Pages2}

    \refer{Paper}{Borgna2015}
	\Rauthor{Borgna J. P., De Leo M., Rial D. and Sanchez de la Vega C.}
	\Rjournal{Commun. Math. Sci, Int. Press Boston, Inc.}
	\Rpages{83-101}
	\Rtitle{General Splitting methods for abstract semilinear evolution equations}
	\Rvolume{13}
	\Ryear{2015}

	\refer{Book}{Cazenave1998}
	\Rauthor{Cazenave T. and Haraux A.}
	\Rpublisher{Oxford Lecture Ser. Math. Appl., Clarendon Press, Rev. ed.}
	\Rtitle{An Introduction to Semilinear Evolution Equations}
	\Ryear{1999}

	\refer{Paper}{DeLeo2015}
	\Rauthor{De Leo M., Rial D. and Sanchez de la Vega C.}
	\Rjournal{IMA J. Numer. Anal.}
	\Rpages{1842-1866}
	\Rtitle{High-order time-splitting methods for irreversible equations}
	\Ryear{2015}

	\refer{Book}{Engel2001}
  \Rtitle{One-parameter semigroups for linear evolution equations}
  \Rauthor{Engel K. J. and Nagel R.}
  \Rvolume{63}
  \Rnumber{2}
  \Rpages{278--280}
  \Ryear{2001}
  \Rpublisher{Springer}

 \refer{Paper}{Fisher1937}
	\Rauthor{Fisher R. A.}
	\Rjournal{Ann. Eugen.}
	\Rpages{353-369}
	\Rtitle{The wave of	advance of advantage	genes}
	\Rvolume{7}
	\Rnumber{4}
	\Ryear{1937}

	\refer{Paper}{Gao2018}
  \Rtitle{Recurrent solutions of the linearly coupled complex cubic-quintic Ginzburg-Landau equations}
  \Rauthor{Gao, P.}
  \Rjournal{Mathematical Methods in the Applied Sciences}
  \Rvolume{41}
  \Rnumber{7}
  \Rpages{2769-2794}
  \Ryear{2018}
  \Rpublisher{Wiley Online Library}

	\refer{Paper}{Ginibre1996}
  \Rtitle{The Cauchy problem in local spaces for the complex Ginzburg-Landau equation I. Compactness methods}
  \Rauthor{Ginibre J. and Velo G.}
  \Rjournal{Physica D: Nonlinear Phenomena}
  \Rvolume{95}
  \Rnumber{3-4}
  \Rpages{191-228}
  \Ryear{1996}
  \Rpublisher{Elsevier}

 \refer{Paper}{Ginibre1997}
  \Rtitle{The Cauchy Problem in Local Spaces for the Complex Ginzburg-Landau Equation II. Contraction Methods}
  \Rauthor{Ginibre J. and Velo G.}
  \Rjournal{Communications in mathematical physics}
  \Rvolume{187}
  \Rnumber{1}
  \Rpages{45-79}
  \Ryear{1997}
  \Rpublisher{Springer}

\refer{Paper}{Guo2000}
  \Rtitle{On existence of almost periodic solution of Ginzburg-Landau equation}
  \Rauthor{Guo, B. and Yuan, R.}
  \Rjournal{Communications in Nonlinear Science and Numerical Simulation}
  \Rvolume{5}
  \Rnumber{2}
  \Rpages{74-79}
  \Ryear{2000}
  \Rpublisher{Elsevier}

\refer{Paper}{Kolmogorov1989}
	\Rauthor{Kolmogorov A. and Petrovsy I. and Piskounov N.}
	\Rjournal{Moscow Univ. Math. Bull.}
	\Rpages{105}
	\Rtitle{Study of the diffusion equation with growth of the quantity of matter and its applications	to a biological problem}
	\Ryear{1989}

		\refer{Paper}{Stepanoff1925}
  \Rtitle{Sur quelques généralisations des fonctions presque périodiques}
  \Rauthor{Stepanoff W.}
  \Rjournal{Comptes rendus}
  \Rvolume{181}
  \Rnumber{20}
  \Rpages{90-92}
  \Ryear{1925}

\refer{Paper}{Weyl1927}
  \Rtitle{Integralgleichungen und fastperiodische Funktionen}
  \Rauthor{Weyl H.}
  \Rjournal{Mathematische Annalen}
  \Rvolume{97}
  \Rnumber{1}
  \Rpages{338-356}
  \Ryear{1927}
  \Rpublisher{Springer}

\end{references}
